\documentclass[11pt,reqno]{amsart}
 \usepackage{setspace}
\usepackage{amsmath,amssymb,amsfonts,url,mathptmx}

\usepackage{color}
\numberwithin{equation}{section}

\newtheorem{thm}{Theorem}[section]

\newcommand{\hdot}{^\text{\r{}}\hspace{-.33cm}H}
\usepackage{lipsum}

\newcommand\blfootnote[1]{%
  \begingroup
  \renewcommand\thefootnote{}\footnote{#1}%
  \addtocounter{footnote}{-1}%
  \endgroup
}

\begin{document}
\title[2$\times$2 singular Liouville Systems]{Degree Counting Theorems for 2$\times$2 non-symmetric singular Liouville Systems } \subjclass{}


\author{YI GU}
\address{Department of Mathematics, University of Florida, Gainesville, FL 32611}

\email{yigu57@ufl.edu}

\begin{abstract} Let $(M,g)$ be a compact Riemann surface with no boundary and $u=(u_1,u_2)$ be a solution of the following singular Liouville system:
$$\Delta_g u_i+\sum_{j=1}^2 a_{ij}\rho_j(\frac{h_je^{u_j}}{\int_M h_je^{u_j}dV_g}-1)=\sum_{l=1}^{N}4\pi\gamma_l(\delta_{p_l}-1), $$
where $h_1,h_2$ are positive smooth functions, $p_1,\cdots,p_N$ are distinct points on $M$, $\delta_{p_l}$ are Dirac masses, $\rho=(\rho_1,\rho_2)(\rho_i\geq 0)$ and $(\gamma_1,\cdots,\gamma_N)(\gamma_l > -1)$ are constant vectors. In the previous work, we derive a degree counting formula for the singular Liouville system when $A$ satisfies standard assumptions. In this article, we establish a more general degree counting formula for 2$\times$2 singular Liouville system when the coefficient matrix $A$ is non-symmetric and non-invertible. Finally, the existence of solution can be proved by the degree counting formula which depends only on the topology of the domain and the location of $\rho$. 

\end{abstract}


\maketitle

\section{Introduction}

In this article we study the following Liouville system defined on a compact Riemann surface $(M,g)$ with no boundary:

\begin{equation}\label{1.1}
\begin{aligned} 
\Delta_g u_1^*+a_{11}\rho_1(\frac{h_1^*e^{u_1^*}}{\int_M h_1^*e^{u_1^*}}-1)+a_{12}\rho_2(\frac{h_2^*e^{u_2^*}}{\int_M h_2^*e^{u_2^*}}-1)=\sum_{l=1}^{N}4\pi\gamma_l(\delta_{p_l}-1), \\
\Delta_g u_2^*+a_{21}\rho_1(\frac{h_1^*e^{u_1^*}}{\int_M h_1^*e^{u_1^*}}-1)+a_{22}\rho_2(\frac{h_2^*e^{u_2^*}}{\int_M h_2^*e^{u_2^*}}-1)=\sum_{l=1}^{N}4\pi\gamma_l(\delta_{p_l}-1),
\end{aligned} 
\end{equation}
where $h_1^*,h_2^*$ are positive smooth functions on $M$, $p_1,\cdots,p_N$ are distinct points on $M$, $\delta_{p_l}$ are Dirac masses, $\rho=(\rho_1,\rho_2)(\rho_i\geq 0)$ and $(\gamma_1,\cdots,\gamma_N)(\gamma_l > -1)$ are constant vectors. Without loss of generality, we assume $Vol(M)=1$. $\Delta_g$ is the Laplace-Beltrami operator ($-\Delta_g\geq 0$). Equation (1.1) is called Liouville system if all the entries in the coefficient matrix $A=(a_{ij})_{2\times2}$ are nonnegative.

In mathematics, physics and many other fields, the Liouville system (\ref{1.1}) has many significant applications. In physics, Liouville systems can be derived form the mean field limit of point vortices of the Euler flow (see \cite{cag-1,cag-2,kie-1,kie-2}). In the classical gauge field theory, Liouville systems are closely related to the Chern-Simons Higgs equation for the non-abelian case\cite{dun,hong,jack,yang}. Here, both the Liouville and Toda equations, together with their solutions, appear prominently in the analysis of the nonrelativistic self-dual Chern-Simons models. Chemotaxis is a prominent feature in the organization of many biological populations. Various Liouville systems are also used to describe such models in the theories of Chemotaxis \cite{chil,kel}.
In geometry, the single Liouville equation is closely related to the famous Nirenberg equation. And the Liouville equation with singular sources describes metrics with conic singularity \cite{kuo}.
In recent developments of related Liouville systems, Chen and Lin completed the program in a series of pioneering works \cite{chenlin1,chenlin2} for Liouville equation. Chipot, Shafrir and Wolansky showed solutions of Liouville systems are radially symmetric with respect to a common point \cite{csw}. Chen-Lin's work was extended by Lin and Zhang to Liouville systems {\cite{linzhang1,linzhang2,linzhang3}}. Then author and Zhang \cite{lei-y} extended Lin-Zhang's work to the systems with Dirac poles. Since our goal is to derive more general degree counting theorem, the main purpose of this article is to extend the existing results (the a priori estimate and degree counting theorem) \cite{lei-y} to 2$\times$2 singular systems with non-symmetric and non-invertible coefficient matrix $A$.

For the coefficient matrix $A$, throughout the paper, we postulate the following conditions:

$(H): a_{11},a_{22}\geq0$, $a_{12},a_{21}>0$ and $a_{21}\geq a_{11},a_{12}\geq a_{22}$.

Obviously, if $u=(u_1,u_2)$ is a solution of (\ref{1.1}), then after adding a constant, $u+c=(u_1+c_1, u_2+c_2)$ is also a solution of (\ref{1.1}). Hence, we can always assume that each component of $u=(u_1,u_2)$ is in 
$$ \hdot^{1}(M):=\{v\in L^2(M);\quad \nabla v\in L^2(M), \mbox{and }\,\, \int_M v dV_g=0\}. $$
Then the equation (\ref{1.1}) is the Euler-Lagrange equation for the following nonlinear functional $J_\rho(u)$ in\hspace{0.1cm} $\hdot^{1}(M)$:
$$J_\rho(u)=\frac{1}{2}\int_M \sum_{i,j=1}^{2} a^{ij}\nabla_g u_i\nabla_g u_j dV_g- \sum_{i=1}^{2}\rho_i \log \int_M h_i e^{u_i} dV_g.$$
Let $\mathbb{N}^+$ be the set of positive integers. To state our theorem, we shall use the following notation:
$$\Sigma:=\{8m\pi+\sum_{p_l\in \Lambda}8\pi(1+\gamma_l);  \Lambda\subset\{p_1,\cdots,p_N\}, m\in \mathbb{N}^+\cup\{0\}\}\setminus \{0\}.  $$
Writing $\Sigma$ as
$$\Sigma=\{8\pi n_k|n_1<n_2<\cdots\},$$
we first establish the following a priori estimate:
\begin{thm}
\label{a priori estimate}
Let $A$ satisfies $(H)$. For $k\in \mathbb{N}^+ \cup \{0\},$ and 
$$\mathcal{O}_k=\{(\rho_1,\rho_2)|\rho_i\geq0,i=1,2; {\rm and}$$
$$8\pi n_k(\frac{a_{21}}{a_{12}}\rho_1+\rho_2)<\frac{a_{11}a_{21}}{a_{12}}\rho_1^2+2a_{21}\rho_1\rho_2+a_{22}\rho_2^2<8\pi n_{k+1}(\frac{a_{21}}{a_{12}}\rho_1+\rho_2) \}.$$
Suppose $h_1,h_2$ are positive and $C^3$ functions on $M$ and $K$ is a compact subset of $\mathcal{O}_k$. Then there exists a constant $C$ such that for any solution $u=(u_1,u_2)$ of (\ref{1.1}) with $\rho \in K$ and $u_i \in \hdot^{1}(M)$, we have
$$|u_i(x)|\le C, \quad \forall x\in M, \quad i=1,2, $$
where $C$ is depended on $M,g,k,K,A,h$. 
\end{thm}

Note that the set $\mathcal{O}_k$ is bounded if all $a_{ii}>0$ and is unbounded if $a_{ii}=0$ for some $i$. By Theorem 1.1, the critical parameter set for (\ref{1.1}) is 
$$\Gamma_k=\{\rho;\frac{a_{11}a_{21}}{a_{12}}\rho_1^2+2a_{21}\rho_1\rho_2+a_{22}\rho_2^2=8\pi n_{k}(\frac{a_{21}}{a_{12}}\rho_1+\rho_2) \}. $$
Thanks to the Chen-Lin's work \cite{chenlin1,chenlin2}, we can define the nonlinear map $T_\rho=(T^1,T^2)$ from $\hspace{0.1cm}  \hdot^{1,2}=\hspace{0.1cm}\hdot^{1}(M)\times\hspace{0.1cm}\hdot^{1}(M)$ to $\hspace{0.1cm}\hdot^{1,2}$ by
$$T^i=-\Delta_g^{-1}(\sum_{j=1}^2 a_{ij}\rho_j(\frac{h_j e^{u_j}}{\int_M h_j e^{u_j}}-1)),\quad i=1,2.$$

Obviously, $T_\rho$ is compact from $\hspace{0.1cm}\hdot^{1,2}$ to itself. Then thanks to the a Priori estimate, for $\rho\notin \Gamma_k$, the Leray-Schauder degree of (\ref{1.1}) can be defined by
$$d_\rho=\deg(I-T_\rho; B_R,0),$$
where $R$ is sufficiently large and $B_R=\{u; u\in\hspace{0.1cm}\hdot^{1,2},\hspace{0.1cm}{
\rm and}\hspace{0.2cm}\sum_{i=1}^2 \|u_i\|_{H^1}<R\}.$ By the homotopic invariance and Theorem 1.1, $d_\rho$ is constant for $\rho \in\mathcal{O}_k$ and is independent of $h=(h_1,h_2)$.

To state the degree counting formula for $d_\rho$, we consider the following generating function $g$:
$$g(x)=(1+x+x^2+\cdots)^{-\mathcal{X}(M)+N}\prod_{l=1}^N (1-x^{1+\gamma_l}), $$
where $\mathcal{X}(M)=2-2g_e(M)$ is the Euler Characteristic of $M$ and  $g_e(M)$ is the genus of $M$. We note that if  $-\mathcal{X}(M)+N< 0$,
$$(1+x+x^2+\cdots)^{-\mathcal{X}(M)+N}=(1-x)^{\mathcal{X}(M)-N}.$$
Writing $g(x)$ in the following form
$$g(x)=1+b_1 x^{n_1}+b_2 x^{n_2}+\cdots+b_k x^{n_k}+\cdots,$$
we use $b_j (j=1,2,\cdots)$ to describe the degree counting theorem:

\begin{thm}
\label{degree counting theorem}
Let $d_\rho$ be the Leray-Schauder degree for (\ref{1.1}). Suppose 
$$ 8\pi n_k(\frac{a_{21}}{a_{12}}\rho_1+\rho_2)<\frac{a_{11}a_{21}}{a_{12}}\rho_1^2+2a_{21}\rho_1\rho_2+a_{22}\rho_2^2<8\pi n_{k+1}(\frac{a_{21}}{a_{12}}\rho_1+\rho_2),$$
then
$$d_\rho=\sum_{j=0}^{k}b_j, \hspace{0.2cm} {\rm where}\hspace{0.2cm} b_0=1.$$
\end{thm}

For most applications $\gamma_l$ are positive integers, which implies that 
$$ \Sigma=\{8\pi m; \hspace{0.1cm}m\in\mathbb{N}^+ \}.$$
Thus in this case if $\mathcal{X}(M) \leq0$ we have
$$g(x)=(1+x+x^2+\cdots)^{-\mathcal{X}(M)}\prod_{l=1}^{N}\frac{1-x^{1+\gamma_l}}{1-x}$$
$$=(1+x+x^2+\cdots)^{-\mathcal{X}(M)}\prod_{l=1}^{N}(1+x+\cdots+x^{\gamma_l})$$
$$=1+b_1x+b_2x^2+\cdots+b_kx^k+\cdots.$$
Obviously $b_j\geq0$ for all $j\geq 1$, which implies
$$d_\rho=1+\sum_{j=1}^{k}b_j>0.$$

\noindent{\bf Corollary 1.1.}
{\em Suppose all $\gamma_l\in \mathbb{N}^+$ and $\mathcal{X}(M) \leq 0$. Then $d_\rho>0$ if 
\begin{equation*}
\frac{a_{11}a_{21}}{a_{12}}\rho_1^2+2a_{21}\rho_1\rho_2+a_{22}\rho_2^2\neq8\pi m(\frac{a_{21}}{a_{12}}\rho_1+\rho_2), \quad \forall m\in \mathbb{N}^+.
\end{equation*}
Thus (\ref{1.1}) always has a solution in this case. }

\medskip

The organization of this article is as follows: In section two we first prove the a priori estimate for the singular Liouville system (\ref{1.1}). When $A$ is non-invertible, we find $u_1$ is a scale multiply of $u_2$. So we can prove the theorem by reducing the system to a single Liouville equation. When $A$ is non-symmetric, we prove the result by reducing it to the symmetric case. That is, $u$ is a solution of singular Liouville system corresponding to the non-symmetric $A$ then we can find a $\tilde u$ which is a solution of the same equation corresponding to a symmetric $\tilde A$. Then we can get the generalization of the a priori estimate and degree counting theorems for this case with some modifications. In section three we prove the degree counting theorem by reducing the system to single Liouville equation and use the previous results of Chen-Lin \cite{chenlin1,chenlin2}. Finally, in section four we discuss some applications of the degree counting theorem for some related topics.

\section{{Proof} of the a priori estimate }
The main purpose of this section is to prove the Theorem 1.1.

Let $u^*=(u_1^*,u_2^*)$ be a solution of (\ref{1.1}). We set 
$$v_i^*=u_i^*-\log\int_M h_i^* e^{u_i^*}dV_g, \hspace{0.1cm} i=1,2.$$
which gives 
$$\int_M h_i^*e^{v_i^*}dV_g=1, \hspace{0.1cm} i=1,2.$$
Therefore, for convenience, in this article we assume that 
\begin{equation}\label{scaleto1}
\int_M h_i^*e^{u_i^*}dV_g=1, \hspace{0.1cm} i=1,2.
\end{equation}
Thus (\ref{1.1}) can be written as
\begin{equation}\label{1.2}
\begin{aligned} 
\Delta_g u_1^*+a_{11}\rho_1(h_1^*e^{u_1^*}-1)+a_{12}\rho_2(h_2^*e^{u_2^*}-1)=\sum_{l=1}^{N}4\pi\gamma_l(\delta_{p_l}-1),\\
\Delta_g u_2^*+a_{21}\rho_1(h_1^*e^{u_1^*}-1)+a_{22}\rho_2(h_2^*e^{u_2^*}-1)=\sum_{l=1}^{N}4\pi\gamma_l(\delta_{p_l}-1).
\end{aligned} 
\end{equation}

Around each singular source, the leading term of $u_i^*$ is a logarithmic function that comes from the following Green's function $G(x,q):$

$$-\Delta_x G(x,q)=\delta_q-1 \quad {\rm and}\quad \int_M G(x,q)dx=0.$$
Define

$$u_i=u_i^*-4\pi \sum_{l=1}^{N} \gamma_l G(x,p_l),$$
and rewrite (\ref{1.2}) as

\begin{equation}\label{manifold}
\begin{aligned} 
\Delta_g u_1+a_{11}\rho_1(h_1e^{u_1}-1)+a_{12}\rho_2(h_2e^{u_2}-1)=0,\\
\Delta_g u_2+a_{21}\rho_1(h_1e^{u_1}-1)+a_{22}\rho_2(h_2e^{u_2}-1)=0,
\end{aligned} 
\end{equation}
where
$$h_i(x)=h_i^*(x)exp\{-\sum_{l=1}^{N}4\pi\gamma_l G(x,p_l)\},$$ 
which implies that around each singular source, say, $p_l$, in local coordinates, $h_j$ can be written as 
$$h_j(x)=|x|^{2\gamma_l}g_j(x)$$
for some positive, smooth function $g_j(x)$.

Let $u=(u_1,u_2)$ be a solution of (\ref{manifold}). To prove a priori estimate for $u$. We only need to prove the upper bound for $u$, because the lower bound of $u$ can be obtained from the upper bound of $u$ and standard Harnack inequality.

Therefore our goal is to prove
\begin{equation}\label{upper bound}
u_i(x)\leq C, \quad i=1,2.
\end{equation}

The proof of (\ref{upper bound}) is by contradiction. Suppose there exists a sequence $u^k$ to (\ref{manifold}) such that $\lim_{k\rightarrow \infty} \max_x\max_i u_i^k(x)\rightarrow \infty$.

Then the equation for $u^k$ is 

\begin{equation}\label{uk}
\begin{aligned} 
\Delta_g u_1^k+a_{11}\rho_1^k(h_1e^{u_1^k}-1)+a_{12}\rho_2^k(h_2e^{u_2^k}-1)=0, \\
\Delta_g u_2^k+a_{21}\rho_1^k(h_1e^{u_1^k}-1)+a_{22}\rho_2^k(h_2e^{u_2^k}-1)=0.
\end{aligned} 
\end{equation}
By an argument similar to a Brezis-Merle type lemma \cite{brez}, it's easy to see that there are only finite blow up points:$\{p_1,\cdots,p_N\}$. And $u_i^k$ is uniformly bounded above in any compact subset away from the blowup set.

Then we consider two cases for the coefficient matrix $A=(a_{ij})$.

{\em Case one:} $det(A)=0$.

$det(A)=0$. This implies $a_{11}a_{22}=a_{12}a_{21}$ and $a_{ij}>0$.

If we multiply $a_{21}$ by the first equation and multiply $a_{11}$ by the second equation in (\ref{uk}) and then subtract, we get
 $$\Delta_g (a_{21} u_1^k-a_{11} u_2^k)=0.$$
Therefore, 

$$u_1^k=\frac{a_{11}}{a_{21}}u_2^k+C^k.$$

By assumption $(H)$, we have either $a_{11}<a_{12}$ or $a_{11}=a_{12}$. Then we are going to discuss this two separate cases.

Case (i): $a_{11}< a_{21}.$

 Let $a=\frac{a_{11}}{a_{21}}(<1)$. Then the second equation in (\ref{uk}) becomes 
\begin{equation}\label{single}
\Delta_g u_2^k+a_{21}\rho_1^k(h_1e^{a u_2^k+C^k}-1)+a_{22}\rho_2^k(h_2e^{u_2^k}-1)=0.
\end{equation}
To apply the local estimate, we rewrite (\ref{single}) in local coordinates. For $p\in M$, let $y=(y^1,y^2)$ be the isothermal coordinates near $p$ such that $y_p(p)=(0,0)$ and $y_p$ depends smoothly on $p$. In this coordinates, $ds^2$ has the form

$$e^{\phi(y_p)}[(dy^1)^2+(dy^2)^2],$$
where 
$$\nabla\phi(0)=0,\phi(0)=0.$$
Also near $p$ we have

$$\Delta_{y_p}\phi=-2Ke^\phi,\quad {\rm where}\hspace{0.1cm} K {\rm \hspace{0.1cm}is\hspace{0.1cm} the \hspace{0.1cm}Gauss\hspace{0.1cm} curvature}. $$
When there is no ambiguity, we write $y=y_p$ for simplicity. In this local coordinates, (\ref{single}) is of the form

$$ -\Delta u_2^k=e^\phi a_{21}\rho_1^k(h_1e^{au_2^k+C^k}-1)+e^\phi a_{22}\rho_2^k(h_2e^{u_2^k}-1),\quad {\rm in}\hspace{0.1cm}B(0,\delta).$$
In this article, we always use $B(p,\delta)$ to denote the ball centered at $p$ with radius $\delta>0$.

Let $f_2^k$ be defined as 

$$ -\Delta f_2^k=-e^\phi a_{21}\rho_1^k-e^\phi a_{22}\rho_2^k,\quad {\rm in}\hspace{0.1cm}B(0,\delta).$$
and $f_2^k(0)=|\nabla f_2^k(0)|=0$. Let $\tilde u_2^k=u_2^k-f_2^k$ and

$$ H_1^k=e^\phi \rho_1^k h_1 e^{af_2^k}, \hspace{0.1cm}H_2^k=e^\phi \rho_2^k h_2 e^{f_2^k},$$
then the equation for $\tilde u_2^k$ becomes

\begin{equation}\label{local1}
-\Delta \tilde u_2^k=a_{21}H_1^ke^{a\tilde u_2^k+C^k}+a_{22}H_2^k e^{\tilde u_2^k},\quad {\rm in} \quad B(0,\delta).
\end{equation}
Set $M_k=\max_x\tilde u_2^k(x)/(1+\gamma)$.
Let
\begin{equation*}
v_2^k(y)=\tilde u_2^k(\epsilon_k y)+2\log\epsilon_k,\quad {\rm where}\hspace{0.1cm} \epsilon_k=e^{-\frac{1}{2}M_k}.
\end{equation*}
Then it is easy to verify that

\begin{equation}\label{rescale}
 -\Delta v_2^k=\epsilon_k^{2-2a} a_{21}H_1^k e^{av_2^k+C^k}+a_{22}H_2^k e^{v_2^k},\quad{\rm in}\quad B(0,\delta\epsilon_k^{-1}).
\end{equation}
Since $u_i^k$ tends to $-\infty$ in $M\backslash \cup_{j=1}^N B(p_j,\delta)$, we have

$$\int_{M\backslash\cup_{j=1}^N B(p_j,\delta)} h_i e^{u_i^k} dV_g \rightarrow 0,\quad i=1,2,$$
and

\begin{equation*}
\begin{aligned} 
 \lim_{k\rightarrow\infty} \int_{B(p_l,\delta)} \rho_1^k h_1 e^{u_1^k} dV_g=\lim_{k\rightarrow\infty} \int_{B(p_l,\delta)} H_1^k e^{a\tilde u_2^k+C^k}dx \\
=\lim_{k\rightarrow\infty}\int_{B(p_l\epsilon_k^{-1},\delta\epsilon_k^{-1})}\epsilon_k^{2-2a}H_1^k e^{av_2^k+C^k}dx=0.
\end{aligned} 
\end{equation*}
Therefore,

$$\int_M h_1 e^{u_1}dV_g = o(1)\neq 1.$$
Contradiction. Thus, blowup cannot happen when $a_{11}<a_{21}$.

\bigskip

Case (ii): $a_{11}=a_{21}$.

In this case, the equation (\ref{uk}) becomes

\begin{equation}
\begin{aligned} \label{uk2}
\Delta_g u_1^k+a_{11}\rho_1^k(h_1e^{u_1^k}-1)+a_{12}\rho_2^k(h_2e^{u_2^k}-1)=0, \\
\Delta_g u_2^k+a_{11}\rho_1^k(h_1e^{u_1^k}-1)+a_{12}\rho_2^k(h_2e^{u_2^k}-1)=0.
\end{aligned} 
\end{equation}
If we subtract this two equation in (\ref{uk2}), we get
$$\Delta_g (u_1^k-u_2^k)=0,$$
which implies

$$ u_1^k= u_2^k+C^k.$$
Thus, we can rewrite  (\ref{uk2}) as a single equation
 
 \begin{equation}\label{single2}
-\Delta_g u_1^k=a_{11}\rho_1^k(h_1e^{u_1^k}-1)+a_{12}\rho_2^k(h_2e^{u_1^k-C^k}-1).
\end{equation}
 (\ref{single2}) in local coordinates can be written as

$$ -\Delta u_1^k=e^\phi a_{11}\rho_1^k(h_1e^{u_1^k}-1)+e^\phi a_{12}\rho_2^k(h_2e^{u_1^k-C^k}-1),\quad {\rm in}\hspace{0.1cm}B(0,\delta).$$
Let $f_1^k$ be defined as 

$$ -\Delta f_1^k=-e^\phi a_{11}\rho_1^k-e^\phi a_{12}\rho_2^k,$$
and $f_1^k(0)=|\nabla f_1^k(0)|=0$. Let $\tilde u_1^k=u_1^k-f_1^k$ and

$$ H_1^k=e^\phi \rho_1^k h_1 e^{f_1^k}, H_2^k=e^\phi \rho_2^k h_2 e^{f_1^k},$$
then the equation for $\tilde u_1^k$ becomes

\begin{equation}
\begin{aligned} 
-\Delta \tilde u_1^k=a_{11} H_1^k e^{\tilde u_1^k}+a_{12} H_2^ke^{\tilde u_1^k-C^k},\quad {\rm in}\hspace{0.1cm}B(0,\delta).
\end{aligned} 
\end{equation}

Again, by an argument similar to a Brezis-Merle type lemma \cite{brez}, there are only finite blow up points:$\{p_1,\cdots,p_N\}$.
Bartolucci, Chen, Lin and Tarantello characterize the blowup solutions of the single Liouville equation (See \cite{bclt,chenlin1,chenlin2}) should satisfy either (Let $\Lambda$ denotes the set of points comes from the singular source.)
 
 \begin{equation*}
\lim_{k\rightarrow\infty}\int_{B(p_l,\delta)} a_{11} H_1^k e^{\tilde u_1^k}+a_{12} H_2^ke^{\tilde u_2^k}dx=8\pi (1+\gamma), \hspace{0.1cm}p_l\in \Lambda,
\end{equation*}
or
 \begin{equation*}
\lim_{k\rightarrow\infty}\int_{B(p_m,\delta)} a_{11} H_1^k e^{\tilde u_1^k}+a_{12} H_2^ke^{\tilde u_2^k}dx=8\pi,\hspace{0.1cm}p_m\notin \Lambda.
\end{equation*}
Here we also observe that
$$\int_{B(p_l,\delta)} H_1^k e^{\tilde u_1^k}dx=\int_{B(p_l,\delta)} \rho_1^k h_1 e^{u_1^k} dV_g,\quad l=1,\cdots,N.$$
Therefore by (\ref{scaleto1}) we have

$$a_{11}\rho_1^k+a_{12}\rho_2^k\rightarrow8\pi n_k.$$
Thus if blowup does happen, $(\rho_1,\rho_2)\in \Gamma_k$.
Therefore if $\rho$ is not on critical hyper-surface $\Gamma_k$, then the a priori estimate holds in this case.

\bigskip
 {\em Case two:} $det(A)\neq0$.
 
 To apply the local estimate, we rewrite the equation (\ref{uk}) in local coordinates

$$ -\Delta u_1^k=e^\phi a_{11}\rho_1^k(h_1e^{u_1^k}-1)+e^\phi a_{12}\rho_2^k(h_2e^{u_2^k}-1),\quad{\rm in}\quad B(0,\delta), $$
$$ -\Delta u_2^k=e^\phi a_{21}\rho_1^k(h_1e^{u_1^k}-1)+e^\phi a_{22}\rho_2^k(h_2e^{u_2^k}-1),\quad{\rm in}\quad B(0,\delta).$$
Let $f_i^k$ be defined as 

$$ -\Delta f_1^k=-e^\phi a_{11}\rho_1^k-e^\phi a_{12}\rho_2^k,$$
$$ -\Delta f_2^k=-e^\phi a_{21}\rho_1^k-e^\phi a_{22}\rho_2^k,$$
and $f_i^k(0)=|\nabla f_i^k(0)|=0$. Let $\tilde u_i^k=u_i^k-f_i^k$ and

$$ H_1^k=e^\phi \rho_1^k h_1 e^{f_1^k}, H_2^k=e^\phi \rho_2^k h_2 e^{f_2^k},$$
then the equation for $\tilde u_i^k$ becomes

\begin{equation}\label{local2}
\begin{aligned} 
\Delta\tilde u_1^k+a_{11}H_1^k e^{\tilde u_1^k}+a_{12}H_2^k e^{\tilde u_2^k}=0, \\
\Delta\tilde u_2^k+a_{21}H_1^k e^{\tilde u_1^k}+a_{22}H_2^k e^{\tilde u_2^k}=0.
\end{aligned} 
\end{equation}

Let

$$B=(b_{ij})_{2\times 2},\hspace{0.3cm} b_{11}=a_{11} \frac{a_{12}}{a_{21}},\hspace{0.5cm}b_{12}=b_{21}=a_{12},\hspace{0.5cm} b_{22}=a_{22},$$
and 

$$U_1=\tilde u_1+\log\frac{a_{21}}{a_{12}},\hspace{0.3cm}U_2=\tilde u_2.$$
Then we can rewrite (\ref{local2}) as

\begin{equation}
\begin{aligned} 
\Delta U_1+b_{11}H_1^ke^{U_1^k}+b_{12}H_2^ke^{U_2^k}=0, \\
\Delta U_2+b_{12}H_1^ke^{U_1^k}+b_{22}H_2^ke^{U_2^k}=0,
\end{aligned} 
\end{equation}
where $B=(b_{ij})_{2\times 2}$ is a symmetric matrix.

Here we observe that

$$\int_{B(p_l,\delta)} H_1^k e^{U_1^k}dx=\frac{a_{21}}{a_{12}} \int_{B(p_l,\delta)} H_1^k e^{\tilde u_1^k}dx=\frac{a_{21}}{a_{12}}\int_{B(p_l,\delta)}\rho_1^kh_1e^{u_1^k}dV_g,$$
$$\int_{B(p_l,\delta)} H_2^k e^{U_2^k}dx= \int_{B(p_l,\delta)} H_2^k e^{\tilde u_2^k}dx=\int_{B(p_l,\delta)}\rho_2^kh_2e^{u_2^k}dV_g.$$
Then we invoke Lemma 2.2 and Lemma 2.3 from \cite{lei-y}, it is easy to see

 \begin{equation}
\int_{M\setminus \cup_{j=1}^N B(p_j,\delta)} h_i e^{u_i^k}dV_g\rightarrow 0,
\end{equation}
 and
 \begin{equation}\label{ratio}
\lim_{k\rightarrow\infty}\int_{B(p_l,\delta)}\rho_i^k h_i e^{u_i^k}dV_g/\mu_l=\lim_{k\rightarrow\infty}\int_{B(p_m,\delta)}\rho_i^k h_i e^{u_i^k}dV_g/\mu_m
\end{equation}
for $i=1,2$ and any pair of $l,m$ between 1 and $N$.

Here we use $\mu_l=1+\gamma_l$ to denote the possible strength of the singular source at each $p_l$ and use $(\sigma_{il})$ to denote the energy  around $p_l$:
$$\sigma_{il}=\lim_{k\rightarrow \infty}\frac{1}{2\pi}\int_{B(p_l,\delta)} H_i^k e^{U_i^k}dx,\hspace{0.1cm} i=1,2.$$
for some $\delta>0.$ Then by (\ref{ratio}) we have, for each $i=1,2$,

$$\frac{\sigma_{i1}}{\mu_1}=\frac{\sigma_{i2}}{\mu_2}=\cdots=\frac{\sigma_{iN}}{\mu_N},$$
and 

$$2\pi(\sigma_{11}+\sigma_{12}+\cdots+\sigma_{1N})=\frac{a_{21}}{a_{12}}\rho_1,$$
$$2\pi(\sigma_{21}+\sigma_{22}+\cdots+\sigma_{2N})=\rho_2.$$
Thus
$$\sigma_{1l}=\frac{a_{21}}{a_{12}}\frac{\rho_1\mu_l}{2\pi\sum_{s=1}^N\mu_s},\quad \sigma_{2l}=\frac{\rho_2\mu_l}{2\pi\sum_{s=1}^N\mu_s},\quad l=1,\cdots,N.$$
For each $l$, the Pohozaev identity (see \cite{lei-y,linzhang1}) for $(\sigma_{1l},\sigma_{2l})$ can be written as

$$\sum_{i,j\in\{1,2\}} a_{ij} \sigma_{il}\sigma_{jl}=4\mu_l\sum_{i\in\{1,2\}} \sigma_{il}.$$
Thus if blowup does happen, $(\rho_1,\rho_2)$ satisfies 

 \begin{equation}
\frac{a_{11}a_{21}}{a_{12}} \rho_1^2+2a_{21}\rho_1\rho_2+a_{22}\rho_2^2=8\pi\sum_{s=1}^{N} \mu_s (\frac{a_{11}}{a_{12}}\rho_1+\rho_2).
\end{equation}
Thus if $\rho$ is not on critical hyper-surfaces $\Gamma_k$, the a priori estimate holds in this case.

\section{{Proof} of degree counting theorems}
The main idea of the proof of the degree counting theorem is to reduce the whole system to the single equation.

{\em Case one:} At least one of $a_{ii}>0$. 

We may assume $a_{11}>0$. Thanks to Theorem 1.1, the Leray-Schauder degree of (\ref{manifold}) for $\rho\in \mathcal{O}_k$ is equal to the degree for the following specific system corresponding to $(\rho_1,0)$:
\begin{equation}
\begin{aligned} 
\Delta_g u_1+a_{11}\rho_1(h_1e^{u_1}-1)=0,\\
\Delta_g u_2+a_{21}\rho_1(h_1e^{u_1}-1)=0,
\end{aligned} 
\end{equation}
where $\rho_1$ satisfies 
$$8\pi n_k<a_{11}\rho_1<8\pi n_{k+1}.$$
It is easy to see that $(\rho_1,0)\in \mathcal{O}_k$, using the degree counting formula of Chen-Lin \cite{chenlin2} for the single equation, we obtain the desired formula in this case:
$$d_\rho=\sum_{j=0}^{k}b_j, \hspace{0.2cm} {\rm where}\hspace{0.2cm} b_0=1.$$

{\em Case two:} $a_{ii}=0$ for both $i=1,2$.

Using $a_{12},a_{21}>0$, we reduce the degree counting formula for $\rho\in \mathcal{O}_k$ to the following system:
\begin{equation}
\begin{aligned} 
\Delta_g u_1+a_{12}\rho_2(h_2e^{u_2}-1)=0,\\
\Delta_g u_2+a_{21}\rho_1(h_1e^{u_1}-1)=0,
\end{aligned} 
\end{equation}
where $\rho_1,\rho_2$ satisfy 
\begin{equation*}
8\pi n_k (\frac{a_{21}}{a_{12}}\rho_1+\rho_2)<2a_{21}\rho_1\rho_2<8\pi n_{k+1}(\frac{a_{21}}{a_{12}}\rho_1+\rho_2).
\end{equation*}
It is easy to see that $(\rho_1,\rho_2)\in \mathcal{O}_k$. Now we consider the special case $\rho_2=\frac{a_{21}}{a_{12}}\rho_1, h_1=h_2=h$. In this case, the maximum principle gives $u_1=u_2+C$. Since they both have average equal to 0, we have $u_1=u_2$. Then (3.2) turns out to be a single equation:
$$\Delta_g u_1+a_{21}\rho_1(he^{u_1}-1)=0,$$
where $\rho_1$ satisfies
$$8\pi n_k<a_{21}\rho_1<8\pi n_{k+1}.$$
By applying the degree counting formula of Chen-Lin \cite{chenlin2} for the single equation, the desired formula can also be obtained in this case:
$$d_\rho=\sum_{j=0}^{k}b_j, \hspace{0.2cm} {\rm where}\hspace{0.2cm} b_0=1.$$
This completes the proof of Theorem 1.2.

\section{{Related} topics}

In this section, we will discuss some applications of the degree counting theorem.

For an open, bounded smooth domain $\Omega$ in $\mathbb{R}^2$, we are interested in the following  Dirichlet problem to the Liouville equations:
\begin{equation}\label{dirichlet}
\left\{\begin{array}{ll}
  \Delta u_i^*+\sum_{j=1}^{2} a_{ij}\rho_j\frac{h_j^*e^{u_j^*}}{\int_\Omega h_j^*e^{u_j^*}}=\sum_{l=1}^{N}4\pi\gamma_l\delta_{p_l}\quad {\rm in}\hspace{0.1cm} \Omega,\\
  u_i|_{\partial \Omega}=0, \quad i=1,2,
  \end{array}\right.
\end{equation}
where $h_1^*,h_2^*$ are smooth functions on $\bar{\Omega}$ and $p_1,\cdots,p_N$ are distinct points in the interior of $\Omega$.

We have the following existence result for (\ref{dirichlet}):
\begin{thm}
\label{dirichlet theorem}
 Suppose 
$$ 8\pi n_k(\frac{a_{21}}{a_{12}}\rho_1+\rho_2)<\frac{a_{11}a_{21}}{a_{12}}\rho_1^2+2a_{21}\rho_1\rho_2+a_{22}\rho_2^2<8\pi n_{k+1}(\frac{a_{21}}{a_{12}}\rho_1+\rho_2),$$
then the Leray-Schauder degree $d_\rho$ for (\ref{dirichlet}) is
$$d_\rho=\sum_{j=0}^{k}b_j, \hspace{0.2cm} {\rm where}\hspace{0.2cm} b_0=1,$$
where $\mathcal{X}(\Omega)=1-g_e(\Omega)$ is the Euler Characteristic of $\Omega$ and  $g_e(\Omega)$ is the number of holes inside $\Omega$. In particular, If $\gamma_1,\cdots,\gamma_N \in \mathbb{N}^+$ and $\Omega$ is not simply connected, we have $d_\rho>0$ and (\ref{dirichlet}) always has a solution.
\end{thm}

\noindent{\bf Remark 4.1.}
{\em To prove Theorem 4.1, we need to show $u^k$ never blow up near the boundary $\partial \Omega$. Since all the singular sources are in the interior of $\Omega$, this fact can be proved by a standard moving plane argument. We refer readers to \cite{linzhang2} for more detailed proof. Then the remaining part of the proof of Theorem 4.1 is the same as Theorem 1.2.}

If the Liouville system on $(M,g)$ is written as 
\begin{equation}\label{special}
  \Delta_g u_i^*+\sum_{j=1}^{2} a_{ij}h_j^*e^{u_j^*}=4\pi\sum_{l=1}^{N}\gamma_l\delta_{p_l},\quad i=1,2,
\end{equation}
with the same assumptions on $h_i^*$ and $vol(M)=1$, $A$ satisfies $(H)$ and $A$ is invertible, here (\ref{special}) is a special case of (\ref{1.1}). Integrating (\ref{special}) on both sides, we get
\begin{equation*}
  \sum_{j=1}^{2} a_{ij}\int_M h_j^*e^{u_j^*}=4\pi\sum_{l=1}^{N}\gamma_l,\quad i=1,2,
\end{equation*}
Thus
\begin{equation}
 \int_M h_j^*e^{u_j^*}=4\pi (\sum_{j=1}^{2} a^{ij})(\sum_{l=1}^{N}\gamma_l),\quad i=1,2,
\end{equation}
where $(a^{ij})$ is the inverse of $(a_{ij})$.

Setting 
$$\rho_i= (\sum_{j=1}^{2}a^{ij})(4\pi\sum_{l=1}^{N}\gamma_l), \quad i=1,2,$$
we can write (\ref{special}) as 
$$ \Delta_g u_1^*+\sum_{j=1}^{2} a_{ij}\rho_j(\frac{h_j^*e^{u_j^*}}{\int_M h_j^*e^{u_j^*}}-1)=\sum_{l=1}^{N}4\pi\gamma_l(\delta_{p_l}-1),\quad i=1,2.$$ 
If $M$ is a torus  ($\mathcal{X}(M)=0$) and $\gamma_l\in N^+$, we can compute the Leray-Schauder degree if $\sum_l \gamma_l$ is odd.
\begin{thm}
Suppose $M$ is a torus, $\gamma_l\in N^+$ and $\sum_l \gamma_l$ is odd. Then the Leray-Schauder degree for  (\ref{special})  is $\frac{1}{2}\Pi_{l=1}^{N}(1+\gamma_l)$.
\end{thm}

\noindent{\bf Proof of Theorem 4.2:} 

Since the genus of the torus $M$ is 1, $\mathcal{X}(M)=0$ and the generating function is 
$$g(x)=\Pi_{p=1}^{N}\frac{1-x^{\mu_p}}{1-x}=\Pi_{p=1}^{N}(1+x+x^2+\cdots+x^{\gamma_p})$$
$$=1+b_1x+b_2x^2\cdots+b_kx^k+\cdots+b_mx^m,$$
where $m=\sum_p\gamma_p$. Let
$$\rho_i=(\sum_{j=1}^{2}a^{ij})(4\pi\sum_{p=1}^{N}\gamma_p),$$
it is easy to see that 
$$ 8\pi n_k(\frac{a_{21}}{a_{12}}\rho_1+\rho_2)<\frac{a_{11}a_{21}}{a_{12}}\rho_1^2+2a_{21}\rho_1\rho_2+a_{22}\rho_2^2<8\pi n_{k+1}(\frac{a_{21}}{a_{12}}\rho_1+\rho_2),$$
for $n_k=(m-1)/2$ and $n_{k+1}=(m+1)/2$. Thus the Leray-Schauder degree $d_\rho$ can be computed as 
$$d_\rho=\sum_{l=0}^{(m-1)/2}b_l.$$
Using $b_{m-1}=b_l$ for $l=0,1,\cdots,m$, we further write $d_\rho$ as 
$$ d_\rho=\frac{1}{2}\sum_{l=1}^mb_l=\frac{g(1)}{2}=\frac{\Pi_{p=1}^N(1+\gamma_p)}{2}.$$
Theorem 4.2 is established.

\end{document}